\documentclass[12pt]{article}
\usepackage{amssymb}
\usepackage{amsfonts}
\usepackage{amsmath}
\usepackage{amsthm}
\usepackage{amscd}
\usepackage{eucal}
\usepackage{latexsym}
\theoremstyle{plain}
\newtheorem{thm}{Theorem}

\theoremstyle{remark}

\numberwithin{thm}{section}
\numberwithin{prop}{section}
\begin{document}
\title{Algebraic Aspects of Multiple Zeta Values}
\author{Michael E. Hoffman\\
\small Dept. of Mathematics\\[-0.8ex]
\small U. S. Naval Academy, Annapolis, MD 21402\\[-0.8ex]
\small \texttt{meh@usna.edu}}
\date{
\small MR Classifications:  Primary 11M06, 11B50; Secondary 05E05, 16W30}
\maketitle
%%AUTHOR-DEFINED MACROS MUST PRECEDE ABSTRACT%%%%%%%%%%%%%%%%
\def\al{\alpha}
\def\be{\beta}
\def\de{\delta}
\def\ep{\epsilon}
\def\zt{\zeta}
\def\la{\lambda}
\def\si{\sigma}
\def\Ga{\Gamma}
\def\De{\Delta}
\def\Si{\Sigma}
\def\del{\nabla}
\def\tilde{\widetilde}
\def\<{\langle}
\def\>{\rangle}
\def\Qxy{\mathbf Q\< x,y\>}
\def\Zxy{\mathbf Z\< x,y\>}
\def\RN{\mathbf R^{\mathbf N}}
\def\Sym{\operatorname{Sym}}
\def\QSym{\operatorname{QSym}}
\def\id{\operatorname{id}}
\def\card{\operatorname{card}}
\def\ev{\operatorname{ev}}
\def\h{\operatorname{ht}}
\def\Zp{\mathbf Z/p\mathbf Z}
\def\A{\mathcal A}
\def\Op{\mathcal O}
\def\H{\mathfrak H}
\def\P{\mathfrak P}
\def\sh{
\setlength{\unitlength}{.5 pt}
\begin{picture}(40,20)
\put(10,2){\line(1,0){20}}
\put(10,2){\line(0,1){10}}
\put(20,2){\line(0,1){10}}
\put(30,2){\line(0,1){10}}
\end{picture}}
%%ABSTRACT%%%%%%%%%%%%%%%%%%%%%%%%%%%%%%%%%%%%%%%%%%%%%%%%%%%%%%
\begin{abstract}
Multiple zeta values have been studied by a wide variety of
methods.
In this article we summarize some of the results about them that
can be obtained by an algebraic approach.
This involves ``coding'' the multiple zeta values by
monomials in two noncommuting variables $x$ and $y$.
Multiple zeta values can then be thought of as defining
a map $\zt:\H^0\to\mathbf R$ from a graded rational vector 
space $\H^0$ generated by the ``admissible words'' of
the noncommutative polynomial algebra $\Qxy$.
Now $\H^0$ admits two (commutative) products making $\zt$ a
homomorphism--the shuffle product and the ``harmonic'' product.
The latter makes $\H^0$ a subalgebra of the algebra $\QSym$ of
quasi-symmetric functions.
We also discuss some results about multiple zeta values
that can be stated in terms of derivations and cyclic derivations 
of $\Qxy$, and we define an action of $\QSym$ on $\Qxy$ that
appears useful.
Finally, we apply the algebraic approach to relations of finite
partial sums of multiple zeta value series.
\end{abstract}
%%ARTICLE TEXT BEGINS HERE%%%%%%%%%%%%%%%%%%%%%%%%%%%%%%%%%%%%%%%
\section{Introduction}
\par
The last fifteen years have seen a great deal of work on the 
multiple zeta values (MZVs)
\begin{equation}
\zt(i_1,i_2,\dots,i_k)=\sum_{n_1>n_2>\cdots>n_k\ge 1}\frac1
{n_1^{i_1}n_2^{i_2}\cdots n_k^{i_k}} ,
\end{equation}
where $i_1,i_2,\dots,i_k$ are positive integers.
The case $k=2$ goes back to Euler \cite{E1}, and was revisited by Nielsen 
\cite{N} and Tornheim \cite{T}.
The general case was introduced in \cite{H1} and \cite{Z}.
These quantities have appeared in a surprising variety of contexts,
including knot theory \cite{LM}, quantum field theory \cite{BG,Kr},
and even mirror symmetry \cite{H4}.
\par
Much work on MZVs has focused on discovering and proving identities
about them, particularly those that express MZVs
of ``length'' ($k$ in equation (1)) greater than one in terms of
ordinary (length one) zeta values.  Even in the length-two case,
it appears that there are MZVs that are ``irreducible'' in the
sense that they can't be expressed (polynomially with rational
coefficients) in terms of length one zeta values, e.g., $\zt(6,2)$.
(Of course it isn't known how to prove even that $\zt(3)^2/\zt(2)^3$
is irrational, so we have to say ``appears'':  but everyone since
Euler who has looked for some reduction of $\zt(6,2)$ hasn't found
one.)
\par
Many approaches have been used to obtain MZV identities.  Analytic
techniques are emphasized in the surveys \cite{BB1} and \cite{BBBL2}.
In this article we will focus on algebraic techniques.  It is
evident that sums of form (1) constitute an algebra by simple multiplication
of series:  this was the starting point of \cite{H2}, which formalized
the ``harmonic algebra'' of MZVs.  But, as has become fairly well known 
by now, there are two distinct algebra structures on the set of MZVs,
the harmonic (or ``stuffle'') algebra and the shuffle algebra.
\par
Before proceeding further it is useful to introduce an algebraic
notation for MZVs.  Series of form (1) can be specified by the
composition (finite sequence of positive integers) $(i_1,i_2,\dots,i_k)$; 
to this composition we assign the word $x^{i_1-1}yx^{i_2-1}y\cdots x^{i_k-1}y$
in noncommuting letters $x$ and $y$.  A series of form (1)
converges exactly when $i_1>1$, i.e., when the corresponding
word starts with $x$ and ends with $y$.  We call such words
``admissible'', and we can think of $\zt$ as assigning
a real number to each admissible word.  (It is convenient to
treat the empty word 1 as admissible and set $\zt(1)=1$.)
Note that if $w$ is the word corresponding to a composition
$(i_1,\dots,i_k)$, the weight $i_1+\dots+i_k$ is the total degree
$|w|$ of $w$.  In this case the length $k$ of the composition
is the $y$-degree of $w$; we denote this by $\ell (w)$.  We will
find it convenient to call $|w|-\ell(w)$ (i.e., the $x$-degree of $w$)
the colength of $w$, denoted $c(w)$.
\par
Let $\H$ be the underlying rational vector space of $\Qxy$, and let 
$\H^0$ be the subspace generated by the admissible words.  Then we
think of $\zt$ as a $\mathbf Q$-linear map $\zt:\H^0\to\mathbf R$.
Now $x$ and $y$ are not admissible, but $\H^0$ is a noncommutative
polynomial algebra on the words $v_{p,q}=x^py^q$ for $p,q\ge 1$
(Of course $\zt$ is not a homomorphism for this algebra structure).
We call the length of a word $w\in\H^0$ in terms of the $v_{p,q}$
its height, denoted $\h(w)$.  For example, $\h(xyx^2y^2)=2$.
\par
With this notation, it is easy to state two identities whose proof
motivated much of the early work on MZVs, the sum theorem and the
duality theorem.  (Both appeared in \cite{H1} as conjectures:  the
sum theorem was proved by Granville \cite{G} and independently by
Zagier; the duality theorem was proved via the iterated integral
discussed below--see \cite{Z}--unfortunately without any notice
of the conjecture!)  The sum theorem can be stated as
$$
\sum_{w\in \H^0,\ |w|=n,\ \ell(w)=k} \zt(w)=\zt(n)
$$
for $n\ge 2$.  For the duality theorem, define an antiautomorphism $\tau$
of the noncommutative polynomial ring $\Qxy$ by $\tau(x)=y$ and
$\tau(y)=x$; note that $\tau$ is an involution that exchanges length
and colength, and preserves height.  The duality theorem states that
$$
\zt(w)=\zt(\tau(w))
$$
for admissible words $w$.
\par
Another of ``early'' results on MZVs was the Le-Murakami theorem of
\cite{LM}.  This is the identity
$$
\sum_{w\in\H^0,\ |w|=2n,\ \h(w)=k}(-1)^{\ell(w)}\zt(w)=
(-1)^n\zt((xy)^n)\sum_{j=0}^{n-k}\binom{2n+1}{2j}(2-2j)B_{2j} ,
$$
which they proved by examining the Kontsevich integral of the unknot.
\par
One reason for the efficacy of the ``algebraic'' notation is apparent--it 
corresponds to the expression of MZVs by iterated integrals, as
follows.  Let $w=a_1a_2\cdots a_n$ be the factorization of an admissible
word into $x$'s and $y$'s.  Then it is easy to show that
\begin{equation}
\zt(w)=\int_0^1\int_0^{t_n}\cdots\int_0^{t_2}\frac{dt_1}{A_n(t_1)}
\cdots\frac{dt_{n-1}}{A_2(t_{n-1})}\frac{dt_n}{A_1(t_n)},
\end{equation}
where
$$
A_i(t)=\begin{cases} t,&\text{if $a_i=x$,}\\1-t,&\text{if $a_i=y$.}
\end{cases}
$$
The duality theorem follows immediately from the change of variable
$(t_1,\dots,t_n)\to (1-t_n,\dots,1-t_1)$ in the iterated integral.
In addition, the fact that iterated integrals multiply via shuffle
product (see the Section 2 below) implies the existence of the shuffle 
product structure on the set of MZVs.
\par
But the series multiplication (or ``stuffle product'') can also
be formated in terms of the algebraic notation; this is the
``harmonic algebra'' of \cite{H2}.  The formulation in \cite{H2}
led to the discovery that the harmonic algebra of MZVs is a 
subalgebra of the quasi-symmetric functions.  We discuss this
in detail in Section 3.
\par
Another remarkable success for the algebraic method is
the result of \cite{H1} that I have since (see
\cite{HO}) called the derivation theorem.  
Let $D$ be the derivation of $\Qxy$ with $D(x)=0$ and $D(y)=xy$.  
Then $D$ takes $\H^0$ to itself, as does the derivation $\tau D\tau$.
We can state the derivation theorem as
$$
\zt(D(w))=\zt(D(\tau(w)))
$$
for admissible words $w$.
The proof of this in \cite{H1} is an elementary but messy partial-fractions
argument.  It seems to have nothing to do with iterated integrals, but
the algebraic notation is working some magic here--just compare the
formulation above with the one given as Theorem 5.1 of \cite{H1}:  for
any admissible composition $(i_1,\dots,i_k)$,
$$
\sum_{j=1}^k\zt(i_1,\dots,i_j+1,\dots,i_k)=
\sum_{\stackrel{\scriptstyle 1\le j\le k}{i_j\ge 2}}
\sum_{p=0}^{i_j-2}\zt(i_1,\dots,i_{j-1},i_j-p,
p+1,i_{j+1},\dots,i_k) .
$$
\par
The sum, duality, and derivation theorems are all subsumed in a
remarkable identity proved in 1999 by Ohno \cite{O}.  It can be
stated nicely in the algebraic notation, but to do so will require
some more machinery:  see Section 4 below.  More recently, I
conjectured, and Ohno proved, a somewhat mysterious ``cyclic''
analogue of the derivation theorem \cite{HO}.  As with the derivation
theorem, the statement in the algebraic notation is very simple, but
the proof is a tricky partial-fractions argument.  We discuss this
in Section 5.
\par
The ``magic'' of the algebraic notation seems to extend to the
finite partial sums of the MZVs.  Here the harmonic algebra still
applies, although the shuffle algebra does not.  In Section 6 we
state some results on finite multiple sums, including some mod
$p$ results ($p$ a prime).  The main result of this section appears
to be new.
\section{The Shuffle Algebra}
\par
As above, let $\H$ be the underlying graded rational vector space of
$\Qxy$, with $x$ and $y$ both given degree 1.  We define a multiplication
$\sh$ on $\H$ by requiring that it distribute over the addition, and 
that it satisfy the following axioms:
\begin{itemize}
\item[S1.]
For any word $w$, $1\sh w=w\sh 1=w$;
\item[S2.]
For any words $w_1,w_2$ and $a,b\in\{x,y\}$,
$$
aw_1\sh bw_2 = a(w_1\sh bw_2)+ b(aw_1\sh w_2).
$$
\end{itemize}
Induction on total degree then establishes the following.
\begin{thm} The $\sh$-product is commutative and associative.
\end{thm}
Recall from the previous section that $\tau$ is the anti-automorphism 
of $\Qxy$ the exchanges $x$ and $y$.  Then we have the following fact.
\begin{thm} $\tau$ is an automorphism of $(\H,\sh)$.
\end{thm}
\begin{proof} Since evidently $\tau^2=\id$, it suffices to show
that $\tau$ is a $\sh$-homomorphism.
Using the axioms S1, S2 above and induction on $|w_1w_2|$,
it is straightforward to prove that
$$
w_1a\sh w_2 b=(w_1\sh w_2b)a + (w_1a\sh w_2)b
$$
for any words $w_1,w_2$ and letters $a,b$.  Now suppose inductively
that $\tau(w_1\sh w_2)=\tau(w_1)\sh\tau(w_2)$ for $|w_1w_2|<n$, and
let $w_1,w_2$ be words with $|w_1w_2|=n$.  We can assume both
$w_1$ and $w_2$ are nonempty; write $w_1=w_1'a$ and $w_2=w_2'b$.
Then
\begin{align*}
\tau(w_1\sh w_2)&=\tau((w_1'\sh w_2)a + (w_1\sh w_2')b)\\
&= \tau(a)\tau(w_1'\sh w_2)+\tau(b)\tau(w_1\sh w_2')\\
&= \tau(a)(\tau(w_1')\sh\tau(w_2))+\tau(b)(\tau(w_1)\sh\tau(w_2'))\\
&= \tau(a)\tau(w_1')\sh\tau(b)\tau(w_2') \\
&= \tau(w_1)\sh\tau(w_2) .
\end{align*}
\end{proof}
\par
Now order the words of $\H$ as follows.  For any words
$w_1,w_2,w_3$, set $w_1xw_2 < w_1yw_3$; and if $u,v$ are words
with $v$ nonempty, set $u<uv$.  A nonempty word $w$ is called
Lyndon if it is smaller than any of its nontrivial right factors;
i.e., $w<v$ whenever $w=uv$ and $u\ne 1\ne v$.  From \cite{Ra}
we have the following result.
\begin{thm} As a commutative algebra, $(\H,\sh)$ is freely generated
by the Lyndon words.
\end{thm}
The link between the shuffle algebra and MZVs is given by the iterated
integral representation (2), together with the well-known fact \cite{Re}
that iterated integrals multiply by shuffle product.  We can state this
as follows.
\begin{thm} The map $\zt: (\H^0,\sh)\to\mathbf R$ is a $\tau$-equivariant
homomorphism.
\end{thm}
The shuffle-product structure has been used to prove some MZV identities.
For example, in \cite{BBBL1} it is first established that 
$$
\sum_{r=-n}^n (-1)^r[(xy)^{n-r}\sh (xy)^{n+r}]=4^n(x^2y^2)^n
$$
in $\H$, and then $\zt$ is applied to get
$$
\sum_{r=-n}^n (-1)^r\zt((xy)^{n-r})\zt((xy)^{n+r})=4^n\zt((x^2y^2)^n).
$$
Using the known result
\begin{equation}
\zt((xy)^k)=\frac{\pi^{2k}}{(2k+1)!}
\end{equation}
(for which see the remarks following Theorem 3.5 below), together with some 
arithmetic, one then obtains
the result conjectured by Zagier \cite{Z} several years earlier:
$$
\zt((x^2y^2)^n)=\frac1{2n+1}\zt((xy)^{2n}) .
$$
Other shuffle convolutions are used to prove some instances of the 
``cyclic insertion conjecture'' for MZVs in the same paper, and the
topic has been revisited in \cite{BB2}.
\section{The Harmonic Algebra and Quasi-Symmetric Functions}
\par
We can define another commutative multiplication $*$ on $\H$ by requiring
that it distribute over the addition and that it satisfy the following
axioms:
\begin{itemize}
\item[H1.]
For any word $w$, $1*w=w*1=w$;
\item[H2.]
For any word $w$ and integer $n\ge 1$,
$$
x^n*w=w*x^n=wx^n;
$$
\item[H3.]
For any words $w_1,w_2$ and integers $p,q\ge 0$,
$$
x^pyw_1*x^qyw_2=x^py(w_1*x^qyw_2)+x^qy(x^pyw_1*w_2)+x^{p+q+1}y(w_1*w_2) .
$$
\end{itemize}
Note that axiom (H3) allows the $*$-product of any pair of words
to be computed recursively, since each $*$-product on the right has
fewer factors of $y$ than the $*$-product on the left-hand side.
Induction on $y$-degree establishes the counterpart of Theorem 2.1.
\begin{thm} The $*$-product is commutative and associative.
\end{thm}
We refer to $\H$ together with its commutative multiplication $*$
as the harmonic algebra $(\H,*)$.  Evidently $\tau$ is \it not \rm
an automorphism of $(\H,*)$.  But we do have counterparts of Theorems
2.3 and 2.4, which are proved in \cite{H2}.
\begin{thm} As a commutative algebra, $(\H,*)$ is freely generated
by the Lyndon words.
\end{thm}
\begin{thm} $(\H^0,*)$ is a subalgebra of $(\H,*)$, and $\zt: (\H^0,*)\to
\mathbf R$ is a homomorphism.
\end{thm}
Because the multiplications $*$ and $\sh$ are quite different, Theorems
2.4 and 3.3 imply that $\zt$ has a large kernel.  For example, since
\begin{align*}
xy*xy &= 2(xy)^2+x^3y \\
xy\sh xy &= 2(xy)^2 + 4x^2y^2
\end{align*}
we must have
$$
\zt(x^3y-4x^2y^2)=0 .
$$
In fact, it has been conjectured that all identities of MZVs come
from comparing the two multiplications.  The derivation theorem
can be recovered, since
$$
y\sh w-y*w=\tau D\tau (w)-D(w)
$$
for $w\in\H^0$ (Theorem 4.3 of \cite{HO}).  Zudilin \cite{Zu} states the
conjecture as 
$$
\ker\zt=\{u\sh v-u*v\ |\ u\in\H^1, v\in\H^0\};
$$
for other formulations see \cite{HJPO} and \cite{W}.
\par
Let $\H^1$ be the vector subspace $\mathbf Q1 +\H y$ of $\H$; it
is evidently a subalgebra of $(\H,*)$.  In fact, since $x$ is the only
Lyndon word ending in $x$, it is easy to see that $\H^1$ is the
subalgebra of $(\H,*)$ generated by the Lyndon words other than
$x$.  Note that any word $w\in\H^1$ can be written in terms of
the elements $z_i=x^{i-1}y$, and that the $y$-degree $\ell(w)$
is the length of $w$ when expressed this way.
We can rewrite the inductive rule (H3) for the $*$-product as
\begin{equation}
z_pw_1*z_qw_2=z_p(w_1*z_qw_2)+z_q(z_pw_1*w_2)+z_{p+q}(w_1*w_2).
\end{equation}
\par
Now for each positive integer $n$, define a map $\phi_n:\H^1\to
\mathbf Q[t_1,\dots,t_n]$ (where $|t_i|=1$ for all $i$) as follows.
Let $\phi_n(1)=1$ and 
$$
\phi_n(z_{i_1}z_{i_2}\cdots z_{i_k})=
\sum_{n\ge n_1>n_2>\dots>n_k\ge 1}t_{n_1}^{i_1}t_{n_2}^{i_2}
\cdots t_{n_k}^{i_k}
$$
for words of length $k\le n$, and let $\phi(w)=0$ for words of
length greater than $n$; extend $\phi_n$ linearly to $\H^1$.
Because the rule (4) corresponds to multiplication of series,
$\phi_n$ is a homomorphism, and $\phi_n$ is evidently
injective through degree $n$.  For each $m\ge n$, there is a
restriction map
$$
\rho_{m,n}:\mathbf Q[t_1,\dots,t_m]\to \mathbf Q[t_1,\dots,t_n]
$$
such that
$$
\rho(t_i)=\begin{cases} t_i,& i\le n\\ 0,& i>n . \end{cases}
$$
The inverse limit
$$
\P=\projlim_n \mathbf Q[t_1,\dots,t_n]
$$
is the subalgebra of $\mathbf Q[[t_1,t_2,\dots]]$ consisting of
those formal power series of bounded degree.  Since the maps $\phi_n$
commute with the restriction maps, they define a homomorphism
$\phi:\H^1\to\P$.
\par
Inside $\P$ is the algebra of symmetric functions
$$
\Sym =\projlim_n \mathbf Q[t_1,\dots,t_n]^{\Si_n}
$$
and also the algebra of quasi-symmetric functions (first
described in \cite{Ges}).
We can define the algebra $\QSym$ of quasi-symmetric functions
as follows.
A formal series $p\in\P$ is in $\QSym$ if the coefficient of 
$t_{i_1}^{p_1}\cdots t_{i_k}^{p_k}$ in $p$ is the same as the
coefficient of $t_{j_1}^{p_1}\cdots t_{j_k}^{p_k}$ in $p$ 
whenever $i_1<i_2<\dots<i_k$ and $j_1<j_2<\dots<j_k$.  Evidently
$\Sym\subset\QSym$.  A vector space basis for $\QSym$ is given by
the monomial quasi-symmetric functions 
$$
M_{(p_1,p_2,\dots,p_k)}=\sum_{i_1<i_2<\dots<i_k}t_{i_1}^{p_1}t_{i_2}^{p_2}
\cdots t_{i_k}^{p_k} ,
$$
which are indexed by compositions $(p_1,\dots,p_k)$.  Since evidently
$\phi(z_{i_1}\cdots z_{i_k})=M_{(i_k,\dots,i_1)}$, we have the following
result.
\begin{thm} $\phi$ is an isomorphism of $\H^1$ onto $\QSym$.
\end{thm}
\par
As is well known, the algebra $\Sym$ of symmetric functions is generated
by the elementary symmetric functions $e_i$, as well as by the power-sum
symmetric functions $p_i$ (Note that we are working over $\mathbf Q$).
It is easy to see that $\phi^{-1}(e_i)=z_1^i$ and $\phi^{-1}(p_i)=z_i$.
Let $\Sym^0$ be the subalgebra of the symmetric functions generated by
the power-sum symmetric functions $p_i$ with $i\ge 2$.  Then
$\phi^{-1}(\Sym)\cap\H^0=\phi^{-1}(\Sym^0)$.  Since $\phi$ is a
homomorphism, we have the following result.
\begin{thm} If $a\in\phi^{-1}(\Sym^0)$, then $\zt(a)$ is a sum of
products of values of $\zt(i)$ of the zeta function with $i\ge 2$.
\end{thm}
\par
In fact, the problem of expressing MZVs $\zt(a)$ with $a\in\phi^{-1}(\Sym^0)$
in terms of values of the zeta function is entirely equivalent to writing
particular monomial symmetric functions in terms of power-sum symmetric
functions $p_i$, for which there are well-known algorithms \cite{Mac}.
This includes cases like
$$
\zt(z_i^k)=\zt(i,i,\dots,i)
$$
(Note $i=2$ occurs in equation (3) above), treated by analytical
methods in \cite{BBB}.  For example, since $M_{22}=\frac12(p_2^2-p_4)$
in $\Sym^0$, we have
$$
\zt(2,2)=\frac12(\zt(2)^2-\zt(4))=\frac12\left(\frac{\pi^4}{36}-
\frac{\pi^4}{90}\right)=\frac{\pi^4}{120} .
$$
(For a general proof of equation (3) by this method, see Corollary 2.3
of \cite{H1}.)
\par
Since $y$ is the only Lyndon word that begins with $y$, we can
write $\H^1=\H^0[y]$ (for either the $\sh$ or the $*$ product).
So we can extend $\zt$ to a map $\hat\zt:\H^1\to\mathbf R$ by
defining $\hat\zt(y)$.  Since 
$$
y*y=2y^2+xy\quad\text{and}\quad y\sh y=2y^2,
$$
there is no way to do this consistently for both multiplications,
but if we restrict our attention to the $*$-multiplication it
turns out that $\hat\zt(y)=\gamma$ (Euler's constant) is a happy choice.
If
$$
H(t)=1+yt+(y^2+xy)t^2+(y^3+yxy+xy^2+x^2y)t^3+\cdots
$$
is the generating function for the complete symmetric functions,
then the following result is easy to show (see \cite{H2}).
\begin{thm}$\hat\zt(H(t))=\Ga(1-t).$
\end{thm}
Now one can show (e.g., using differential equations) that
$$
\sum_{w\in\H^0,\ \h(w)=1} \zt(w)u^{c(w)}v^{\ell(w)}=
1-\frac{\Ga(1-u)\Ga(1-v)}{\Ga(1-u-v)} .
$$
Putting this together with Theorem 3.6, we have
\begin{equation}
\sum_{w\in\H^0,\ \h(w)=1} \zt(w)u^{c(w)}v^{\ell(w)}=
\zt\left(1-\frac{H(u)H(v)}{H(u+v)}\right) .
\end{equation}
Hence $\zt(w)\in\zt(\phi^{-1}(\Sym^0))$ for any word $w$ of height 1 (i.e.,
of the form $x^py^q$), and thus can be written in terms of ordinary
zeta values $\zt(n)=\zt(z_n)$.
\par
Remarkably, Ohno and Zagier \cite{OZ} have recently proved that
equation (5) is just the constant term of the following result.
\begin{thm}
$$
\sum_{w\in\H^0}\zt(w)u^{c(w)}v^{\ell(w)}z^{\h(w)-1}=
\frac1{1-z}\zt\left(1-\frac{H(u)H(v)}{H(\al)H(\be)}\right),
$$
where 
\begin{align*}
\al &= \frac12\left((u+v)+\sqrt{(u+v)^2-4uvz}\right) \\
\be &= \frac12\left((u+v)-\sqrt{(u+v)^2-4uvz}\right) .
\end{align*}
\end{thm}
The theorem implies that any sum of MZVs of fixed weight, length, and height,
e.g.,
$$
\sum_{|w|=6,\ \h(w)=2,\ \ell(w)=3}\zt(w)=\zt(3,2,1)+\zt(2,3,1)+\zt(2,1,3)
+\zt(3,1,2)
$$
is in $\zt(\phi^{-1}(\Sym^0))$ and hence expressible in terms of $\zt(n)$'s.
But the theorem implies much more.  For example, taking the limit as
$z\to 1$ gives the sum theorem, and setting $v=-u$ gives the Le-Murakami
theorem.
\par
For another application of Theorem 3.6 see \cite{H4}.
\section{Derivations and an Action by Quasi-Symmetric Functions}
\par
As mentioned in the introduction, the derivation theorem has a
far-reaching generalization proved by Ohno \cite{O}.  In this section
we give a succinct statement of Ohno's theorem and some of its equivalents
using the Hopf algebra structure of $\QSym$.  
\par
We begin by motivating the use of a Hopf algebra structure in this context.
(The standard references on Hopf algebras are \cite{S} and \cite{MM},
but the reader may find a source like \cite{Ka} more convenient.)
Let $\Op$ be an algebra of operators (with composition as multiplication)
acting on an algebra $\A$.  Then elements of the tensor product $\Op\otimes\Op$
act naturally on products $pq$ for $p,q\in\A$:  
$\al\otimes\be(pq)=\al(p)\be(q)$.  
To say that $\al\in\Op$ is a derivation is to say that the action of $\al$
on a products agrees with the action of $\al\otimes 1+1\otimes\al$:
$$
\al(pq)=\al(p)q+p\al(q)=(\al\otimes 1+1\otimes\al)(pq) .
$$
A Hopf algebra structure on $\Op$ is essentially a ``coproduct''
$\De:\Op\to\Op\otimes\Op$ compatible with the multiplication in $\Op$.
We require that $\al(pq)=\De(\al)(pq)$ for all $\al\in\Op$.
Elements $\al$ with $\De(\al)=\al\otimes 1+1\otimes\al$ are called
primitive, so the primitives in $\Op$ are exactly those that act as
derivations.  The ``fine print'' of the definition of a (graded connected) 
Hopf algebra requires that $\De(\al)$ always contain the terms $\al\otimes 1$
and $1\otimes\al$ for $\al$ of positive degree, so primitive elements 
are those whose coproducts are as simple as possible.  We can generalize
the notion of derivation by allowing extra terms in the coproduct.
For example, a set $\{\al_0=1,\al_1,\al_2,\dots\}$ of elements is called a
set of divided powers if
$$
\De(\al_n)=\sum_{i+j=n}\al_i\otimes \al_j;
$$
if we think of the $\al_n$ as operators, they are sometimes called a
``higher derivation''.  Thus, a Hopf algebra of operators is a natural 
extension of the notion of a Lie algebra acting by derivations.
\par
Now $(\H^1,*)\cong\QSym$ has a Hopf algebra structure with coproduct
$\De$ defined by
$$
\De(z_{i_1}z_{i_2}\cdots z_{i_n})=
\sum_{j=0}^n z_{i_1}\cdots z_{i_j}\otimes z_{i_{j+1}}\cdots z_{i_n} ,
$$
(and counit $\ep$ with $\ep(u)=0$ for all elements $u$ of positive degree).
This extends the well-known Hopf algebra structure on the algebra
$\Sym$ (as described in \cite{Gei}), in which the elementary symmetric
functions $e_i$ ($\leftrightarrow y^i$) and complete symmetric functions
$h_i$ are divided powers, while the power sums $p_i$ ($\leftrightarrow z_i$)
are primitive.  The Hopf algebra $(\H^1,*,\De)$ is commutative but not
cocommutative.  Its (graded) dual is the Hopf algebra of noncommutative
symmetric functions as defined in \cite{Gel}.
\par
Now define $\cdot:\H^1\otimes\Qxy\to\Qxy$ by setting $1\cdot w=w$
for all words $w$, 
$$
z_k\cdot 1=0,\quad z_k\cdot x=0,\quad z_k\cdot y=x^ky
$$
for all $k\ge 1$, and 
\begin{equation}
u\cdot w_1w_2=\sum_u (u'\cdot w_1)(u''\cdot w_2)
\end{equation}
where $\De(u)=\sum_u u'\otimes u''$; the coassociativity of $\De$ insures
this is well-defined.  It turns out (Lemma 5.2 of \cite{HO})
that $u\cdot w$ just consists
of those terms of $u*w$ having the same $y$-degree as $w$, 
so it follows (from the associativity of $*$)
that $\cdot$ is really an action, i.e., $u\cdot(v\cdot w)=(u*v)\cdot w$.
Also, equation (6) says the action makes $\Qxy$ a $\QSym$-module
algebra, in the terminology of \cite{Ka}.
\par
We note that the action of $z_1$ on $\Qxy$ is just the derivation 
$D$ defined in the introduction, since $z_1\cdot x=0$ and $z_1\cdot y=xy$.
In fact, for each $n\ge 1$ we have a derivation $D_n$ given by
$D_n(w)=z_n\cdot w$, since the $z_n$ are primitive in $\QSym$.
\par
In terms of this action, we can now state Ohno's theorem \cite{O} as follows.
\begin{thm} For any word $w\in\H^0$ and nonnegative integer $i$,
$$
\zt(h_i\cdot w)=\zt(h_i\cdot\tau(w)).
$$
\end{thm}
Recall that the $h_n$ are divided powers, i.e., 
$\De(h_n)=\sum_{i+j=n}h_i\otimes h_j$.  Only $h_1=z_1$ is primitive,
in which case we recover the derivation theorem.  Taking $h_0=1$ gives
the duality theorem, and with a little manipulation the sum theorem
can also be obtained.
\par
M. Kaneko sought to generalize the derivation theorem in another
way.  One can formulate the derivation theorem as saying that
$(\tau D\tau-D)(w)\in\ker\zt$ for all $w\in\H^0$.  Are there derivations
of higher degree for which this is still true?  Kaneko defined a
degree-$n$ derivation $\partial_n$ of $\Qxy$ by
$$
\partial_n(x)=-\partial_n(y)=x(x+y)^{n-1}y ,
$$
and conjectured that $\partial_n(w)\in\ker\zt$ for all $w\in\H^0$.
Note $\partial_1=\tau D\tau-D$, so the conjecture holds for $n=1$;
and the case $n=2$ follows easily from Theorem 4.1.
\par
Eventually Kaneko and K. Ihara proved the conjecture \cite{IK} by showing 
it equivalent to Theorem 4.1.  One way to see this involves the action
we have just defined.  Extend the action of $\QSym$ on $\H$ to an
action of $\QSym[[t]]$ on $\H[[t]]$ in the obvious way, and (as in
the previous section) let
$$
H(t)=1+h_1t+h_2t^2+\cdots\in\QSym[[t]] 
$$
be the generating function of the complete symmetric functions.
If we set $\si_t(u)=H(t)\cdot u$ for $u\in\H$, then Theorem 4.1 is
equivalent to $\zt(\bar\si_t (u)-\si_t(u))=0$ for $u\in\H$, where
$\bar\si_t=\tau\si_t\tau$.
Now $\si_t$ is an automorphism of $\H^0[[t]]$:  in fact $\si_t^{-1}(u)
=E(-t)\cdot u$, where
$$
E(t)=1+yt+y^2t+\cdots\in\QSym[[t]]
$$
is the generating function of the elementary symmetric functions.
Thus, Theorem 4.1 is equivalent to 
$$
\bar\si_t\si_t^{-1}(u)-u\in\ker\zt
$$
for all $u\in\H^0[[t]]$.  Then following result implies Kaneko's
conjecture.
\begin{thm}
$$
\bar\si_t\si_t^{-1}=\exp\left(\sum_{n=1}^\infty\frac{t^n}{n}\partial_n
\right) .
$$
\end{thm}
This result can be proved by showing both sides are automorphisms of $\H[[t]]$
that fix $t$ and $x+y$, and take $x$ to $x(1-ty)^{-1}$ (see \cite{HO}).
The derivations $\partial_n$ are related to the derivations $D_n$ mentioned
above as follows.  Since
$$
\frac{d}{dt}\log H(t)=\frac{H'(t)}{H(t)}=\sum_{n=1}^\infty p_n t^{n-1},
$$
the map $\si_t$ can also be written
$$
\si_t=\exp\left(\sum_{n=1}^\infty \frac{t^n}{n} D_n\right) .
$$
Hence Theorem 4.2 says that 
$$
\exp\left(\sum_{n=1}^\infty\frac{t^n}{n}\partial_n\right)=
\exp\left(\sum_{n=1}^\infty \frac{t^n}{n}\bar D_n\right) 
\exp\left(-\sum_{n=1}^\infty \frac{t^n}{n}D_n\right) ,
$$
where $\bar D_n=\tau D_n\tau$.  Thus, the $\partial_n$ can be written
in terms of the $D_n$ and $\bar D_n$ via the Campbell-Hausdorff formula.
For example,
$$
\partial_2=\bar D_2-D_2-[\bar D_1,D_1],
$$
and 
$$
\partial_3=\bar D_3-D_3-\frac34[\bar D_1,D_2]-\frac34[\bar D_2,D_1]
+\frac14[[\bar D_1,D_1],D_1]-\frac14[\bar D_1,[\bar D_1,D_1]].
$$
\section{Cyclic Derivations}
\par
There is an analogue of the derivation theorem involving a ``cyclic
derivation'' $C:\H\to\H$.  We can define $C$ as the composition
$\tilde\mu\hat C$, where $\hat C:\H\to\H\otimes\H$ is the derivation
sending $x$ to 0 and $y$ to $y\otimes x$, and $\tilde\mu(a\otimes b)=ba$.
Here we regard $\H\otimes\H$ as a two-sided module over $\H$ via
$a(b\otimes c)=ab\otimes c$ and $(a\otimes b)c=a\otimes bc$.  Thus, e.g.,
\begin{align*}
C(x^3yxy) &= \tilde\mu(x^3(y\otimes x)xy + x^3yx(y\otimes x))\\
&= \tilde\mu(x^3y\otimes x^2y+x^3yxy\otimes x)\\
&= x^2yx^3y+x^4yxy .
\end{align*}
This particular definition follows D. Voiculescu's version of the cyclic
derivative \cite{V}:  cyclic derivatives were first studied by Rota,
Sagan and Stein \cite{RSS}.
\par
In terms of the composition notation, $C$ differs from $D$ in that the
entries are permuted cyclically, e.g.,
$$
D(4,2)=(5,2)+(4,3)\quad\text{versus}\quad C(4,2)=(5,2)+(3,4).
$$
The following result was conjectured by myself and proved by Ohno
\cite{HO}.
\begin{thm} For any word $w\in\H^1$ that is not a power of $y$,
$$
\zt(C(w))=\zt(\tau C\tau(w)).
$$
\end{thm}
As mentioned in the introduction, the proof uses partial fractions.
\par
The difference between $C$ and $D$ is most striking when applied to
periodic words.  For example, Theorem 5.1 applied to $w=(x^2y)^n$ gives
(in the composition notation)
$$
\zt(4,3,\dots,3)=\zt(3,3,\dots,3,1)+\zt(2,3,\dots,3,2).
$$
Theorem 5.1 also gives a very nice proof of the sum theorem.  Here
is the idea:  Let $u=x+ty$.  Then the coefficient of $t^k$ in $xu^{n-2}y$
is the sum of all words $w\in\H^0$ with $|w|=n$ and $\ell(w)=k$.
Now
$$
C(u^{n-1})=(n-1)txu^{n-2}y\quad\text{while}\quad \tau C\tau(u^{n-1})
=(n-1)xu^{n-2}y,
$$
so the cyclic derivation theorem implies that $\zt$ applied to the 
coefficient of $t^{k-1}$ equals $\zt$ applied to the coefficient of
$t^k$. That is, the sum of MZVs of fixed weight $n$ and length $k$ 
must be independent of $k$ (and so must be $\zt(n)$).
\par
Here is another corollary of Theorem 5.1, stated in terms of the
action of $\QSym$ on $\Qxy$.
\begin{thm} For $m,n\ge 1$, $\zt(z_n\cdot xy^m)=\zt(z_m\cdot xy^n)$.
\end{thm}
\par
Just as the derivation theorem extends to Theorem 4.1, it is natural
to ask if the cyclic derivation theorem can be extended.  It is easy
to define cyclic derivations $C_n$ analogous to the $D_n$ of the last
section:  just set $C_n=\tilde\mu\hat C_n$, where $\hat C_n(x)=0$ and
$\hat C_n(y)=y\otimes x^n$.  One could then try to define cyclic
derivations analogous to Kaneko's derivations $\partial_n$ (which
are expressible in terms of commutators of the $D_n$ and $\bar D_n$).
The difficulty appears to be in defining the commutator of cyclic
derivations.
\section{Finite Multiple Sums and Mod $p$ Results}
\par
In this section we consider the finite sums
$$
A_{(i_1,\dots,i_k)}(n)=
\sum_{n\ge n_1>n_2>\dots>n_k\ge 1}\frac1{n_1^{i_1}\cdots n_k^{i_k}}
$$
and
$$
S_{(i_1,\dots,i_k)}(n)=
\sum_{n\ge n_1\ge n_2\ge\dots\ge n_k\ge 1}\frac1{n_1^{i_1}\cdots n_k^{i_k}}
$$
where the notation is patterned after that of \cite{H1}; the multiple
zeta values of the previous sections are
$$
\zt(i_1,\dots,i_k)=\lim_{n\to\infty} A_{(i_1,\dots,i_k)}(n),
$$
when the limit exists (i.e., when $i_1>1$).
\par
The sums $A_I(n)$ and $S_I(n)$ are related in an obvious way, e.g.,
$$
S_{(4,2,1)}(n)=A_{(4,2,1)}(n)+A_{(6,1)}(n)+A_{(4,3)}(n)+A_{(7)}(n) .
$$
We can formalize the relation as follows.  For compositions $I,J$,
we say $I$ refines $J$ (denoted $I\succ J$) if $J$ can be obtained
from $I$ by combining some of its parts.  Then
\begin{equation}
S_I(n)=\sum_{I\succeq J} A_J(n) .
\end{equation}
Of course $S_{(m)}(n)=A_{(m)}(n)$ for all $m,n$.
\par
It will be useful to have some additional notations for compositions.
We adapt the notation used in previous sections for words, so 
for $I=(i_1,\dots,i_k)$ the weight of $I$ is $|I|=i_1+\dots+i_k$,
and $k=\ell(I)$ is the length of $I$.  
For $I=(i_1,\dots,i_k)$, the reversed composition $(i_k,\dots,i_1)$ will
be denoted $\bar I$:  of course reversal preserves weight, length and
refinement (i.e., $I\succeq J$ implies $\bar I\succeq\bar J$).
\par
Compositions of weight $n$ are in 1-to-1 correspondence with subsets of 
$\{1,2,\dots,n-1\}$ via partial sums
$$
(i_1,i_2,\dots,i_k)\to \{i_1,i_1+i_2,\dots,i_1+\dots+i_{k-1}\},
$$
and $I\preceq J$ if and only if the subset corresponding to $I$ contains
that corresponding to $J$.  Complementation in the power set then
gives rise to an involution $I\to I^*$; e.g., $(1,1,2)^*=(3,1)$.
Evidently $|I^*|=|I|$ and $\ell(I)+\ell(I^*)=|I|+1$.  Also, $I\preceq J$
if and only if $I^*\succeq J^*$.  Finally, for two compositions $I$
and $J$ we write $I\sqcup J$ for their juxtaposition.
\par
From \cite{H1} we have formulas for symmetric sums of $A_I(n)$ and
$S_I(n)$ in terms of length one sums $S_{(m)}(n)$.  (Though the proofs
in \cite{H1} are given for infinite series, they carry over to the 
finite case.)  They require some notation to state.  For a partition
$\Pi=\{P_1,\dots,P_l\}$ of the set $\{1,2,\dots,k\}$, let 
$$
c(\Pi)=\prod_{s=1}^l(\card P_s-1)!\quad\text{and}\quad
\tilde c(\Pi)=(-1)^{k-l}\prod_{s=1}^l(\card P_s-1)! ,
$$
and if also $I=(i_1,\dots,i_k)$ is a composition of length $k$, let 
$$
S(n,\Pi,I)=\prod_{s=1}^l S_{(p_s)}(n) ,\quad\text{where}\quad
p_s=\sum_{j\in P_s} i_j .
$$
If $I=(i_1,\dots,i_k)$ is a composition of
length $k$, then elements $\si\in\Si_k$ of the symmetric group
act on $I$ via $\si\cdot I=(i_{\si(1)},\dots,i_{\si(k)})$.
Then Theorems 2.1 and 2.2 of \cite{H1} give us the following result.
\begin{thm} For all positive integers $k$ and $n$ and compositions
$I$ of length $k$,
\begin{align*}
\sum_{\si\in\Si_k}S_{\si\cdot I}(n)&=
\sum_{\text{partitions of $\{1,\dots,k\}$}}c(\Pi)S(n,\Pi,I)\\
\sum_{\si\in\Si_k}A_{\si\cdot I}(n)&=
\sum_{\text{partitions of $\{1,\dots,k\}$}}\tilde c(\Pi)S(n,\Pi,I)
\end{align*}
\end{thm}
\par
Because of the correspondence between compositions and noncommutative
words in $\H^1$, we have (for any fixed $n$) a map $\rho_n:\H^1\to\mathbf Q$
sending $w\in\H^1$ to $A_{I(w)}(n)$, where $I(w)$ is the composition 
associated with $w$.
Note that $\rho_n$ is the composition $\ev\circ T\circ\phi_n$, where
$\phi_n$ is the map defined in Section 3, $T$ is the automorphism
of $\QSym$ sending $M_I$ to $M_{\bar I}$, and $\ev$ is the function
that sends $t_i$ to $\frac1{i}$.  Thus, $\rho_n:(\H^1,*)\to\mathbf R$ 
is a homomorphism.  
We can combine the homomorphisms $\rho_n$ into a homomorphism
$\rho$ that sends $w\in\H^1$ to the real-valued sequence $n\to\rho_n(w)$.
We shall write $A_I$ for the real-valued sequence $n\to A_I(n)$
(and similarly for $S_I$), so $\rho$ sends $w$ to $A_{I(w)}$.
\par
Now $\QSym$ has various integral bases besides the $M_I$.  In the
literature one often sees the fundamental quasi-symmetric functions
$$
F_I=\sum_{J\succeq I} M_J,
$$
but we will be concerned with what we call the ``essential'' quasi-symmetric
functions 
$$
E_I=\sum_{J\preceq I} M_J .
$$
In view of equation (7), the homomorphism $\rho_n$ sends $E_I$ to $S_I(n)$.
\par
Since $\QSym$ is a commutative Hopf algebra, its antipode
$S$ is an automorphism of $\QSym$ and $S^2=\id$.  Now $S$ can
be given by the following explicit formulas:  for proof see
\cite{Eh} or \cite{H3}.
\begin{thm} The antipode $S$ of $\QSym$ is given by 
\begin{itemize}
\item[1.]
$S(M_I)=\sum_{I_1\sqcup I_2\sqcup\cdots\sqcup I_l=I} (-1)^l M_{I_1}M_{I_2}
\cdots M_{I_l};$
\item[2.]
$S(M_I)=(-1)^{\ell(I)}E_{\bar I}$.
\end{itemize}
\end{thm}
Part (2) of this result says that the $E_I$ have essentially the
same multiplication rules as the $M_I$:  if $T$ is the automorphism of 
$\QSym$ defined above, then $S\circ T$ takes  any identity among the $M_I$ 
to an identity among the $E_I$ that differs only in signs.  
For example, since
$$
M_{(2)}M_{(3)}=M_{(2,3)}+M_{(3,2)}+M_{(5)}
$$
we have
$$
E_{(2)}E_{(3)}=E_{(2,3)}+E_{(3,2)}-E_{(5)} .
$$
\par
Now define an automorphism $\psi$ of $\Qxy$ by 
$$
\psi(x)=x+y,\quad \psi(y)=-y
$$
Evidently $\psi^2=\id$, and $\psi(\H^1)=\H^1$.
Thus $\psi$ defines a linear involution of $\H^1\cong
\QSym$ (which is \it not\rm, however, a homomorphism for the 
$*$-product).  We can describe the action of $\psi$ on the integral
bases for $\QSym$ as follows.
\begin{thm} For any composition $I$,
\begin{itemize}
\item [1.]
$\psi(M_I)=(-1)^{\ell(I)}F_I$
\item [2.]
$\psi(E_I)=-E_{I^*}$
\end{itemize}
\end{thm}
\begin{proof} Suppose $w=w(I)$ is the word in $x$ and $y$ corresponding
to a composition $I$.  Then evidently substituting $y$ in place of any
particular factor $x$ in $w$ corresponds to splitting a part of $I$.  
With this observation, part (1) is clear (there is also one factor of 
$-1$ for each occurrence of $y$ in $w$).
\par
Now we prove part (2).  We have
$$
\psi(E_I)=\sum_{J\preceq I}\psi(M_J)=
\sum_{J\preceq I} (-1)^{\ell(J)}F_J
$$
from part (1).  From Example 1 of \cite{H3}, $S(F_I)=(-1)^{|I|}
F_{\bar I^*}$, where $S$ is the antipode of $\QSym$.  Thus
\begin{multline*}
S\psi(E_I)=\sum_{J\preceq I}(1)^{\ell(J)+|J|}F_{\bar J^*}
=-\sum_{J\preceq I}(-1)^{\ell(J^*)}F_{\bar J^*}\\
=-\sum_{\bar J^*\succeq\bar I^*}(-1)^{\ell(\bar J^*)}F_{\bar J^*}
=-\sum_{K\succeq\bar I^*}(-1)^{\ell(K)}F_K .
\end{multline*}
Now by M\"obius inversion,
$$
F_I=\sum_{I\preceq J}M_J\quad\text{implies}\quad
M_I=\sum_{I\preceq J}(-1)^{\ell(I)-\ell(J)}F_J ,
$$
and so
$$
S\psi(E_I)=-(-1)^{\ell(\bar I^*)} M_{\bar I^*}
$$
Apply $S$ be both sides to get
$$
\psi(E_I)=-(-1)^{\ell(I^*)}(-1)^{\ell(\bar I^*)}E_{I^*}=-E_{I^*}.
$$
\end{proof}
\par
We consider two operators on the space $\RN$ of real-valued sequences.
First, there is the partial-sum operator $\Si$, given by
$$
\Si a(n) = \sum_{i=0}^n a(i) 
$$
for $a\in\RN$.
Second, there is the operator $\del$ given by
$$
\del a(n) = \sum_{i=0}^n \binom{n}{i}(-1)^i a(i) .
$$
It is easy to show that $\Si$ and $\del$ generate a dihedral group
within the automorphisms of $\RN$, i.e., $\del^2=\id$ and $\Si\del=
\del\Si^{-1}$.  It follows that $(\Si\del)^2=\id$.  
We have the following result on multiple sums.
\begin{thm} For any composition $I$, $\Si\del S_I=-S_{I^*}$.
\end{thm}
\begin{proof} We proceed by induction on $|I|$.  The weight one
case is $\Si\del S_{(1)}=\del\Si^{-1}S_{(1)}=-S_{(1)}$, i.e.
$$
\sum_{k=1}^n\frac{(-1)^k}{k}\binom{n}{k}=-\sum_{k=1}^n\frac1{k} ,
$$
which is a classical (but often rediscovered) formula; it
actually goes back to Euler \cite{E2}.  For $I=(i_1,i_2,\dots,i_k)$,
it is straightforward to show that $\del S_I(n)=\frac1{n}\del f(n)$, 
where $f\in\RN$ is given by
$$
f(n)=\begin{cases} S_{(i_2,\dots,i_k)}(n),&\text{if $i_1=1$;}\\
\Si^{-1}S_{(i_1-1,i_2,\dots,i_k)}(n),&\text{otherwise.}\end{cases}
$$
\par
Now suppose the theorem has been proved for all $I$ of weight less
than $n$, and let $I=(i_1,\dots,i_k)$ have weight $n$.  There are
two cases:  $i_1=1$, and $i_1>1$.  In the first case, let
$(i_2,\dots,i_k)^*=J=(j_1,\dots,j_r)$.  By the assertion of the
preceding paragraph and the induction hypothesis,
$$
\Si\del S_I(n)=\Si(\frac1{n}\del S_{J^*}(n))=-\Si(\frac1{n}\Si^{-1}S_J(n))
=-S_{(j_1+1,j_2,\dots,j_r)}(n).
$$
But evidently $I^*=(j_1+1,j_2,\dots,j_r)$, so the theorem holds in 
this case
\par
If $i_1>1$, we instead write $(i_1-1,i_2,\dots,i_k)^*=J=(j_1,\dots,j_r)$.
Then
$$
\Si\del S_I(n)=\Si(\frac1{n}\del\Si^{-1}S_{J^*}(n))=\Si(\frac1{n}\Si\del
S_{J^*}(n))=-\Si(\frac1{n}S_J(n))=-S_{(1,j_1,\dots,j_r)}(n).
$$
But in this case $I^*=(1,j_1,\dots,j_r)$, so the theorem holds in this
case as well.
\end{proof}
\par
The proof of the preceding result is essentially a formalization of
the procedure in App. B of \cite{Ve}.  (For a recent occurrence of
the special case $I=(1,1,1)$ as a problem, see \cite{He}.)  
Theorem 6.4, together with part (2) of Theorem 6.3, says that the 
diagram
\begin{equation}
\begin{CD}
\QSym @>{\psi}>>\QSym\\
@V{\rho}VV @V{\rho}VV\\
\RN @>{\Si\del}>>\RN
\end{CD}
\end{equation}
commutes.
\par
For the rest of this section, we discuss mod $p$ results about
$S_I(p-1)$ and $A_I(p-1)$, where $p$ is a prime. (Some results
of this type appear in \cite{Zh}.)  For prime $p$, the sums
$A_I(p-1)$ and $S_I(p-1)$ contain no factors of $p$ in the
denominators, and can be regarded as elements of the field
$\Zp$.  The following result about length one harmonic sums is
well known (cf. \cite{HW}, pp. 86-88).
\begin{thm} $S_{(k)}(p-1)\equiv 0 \mod p$ for all prime $p>k+1$.
\end{thm}
Because Theorem 6.1 expresses symmetric sums of $S_I(p-1)$ and
$A_I(p-1)$ in terms of length one sums, any such symmetric sum
is zero mod $p$ for $p>|I|+1$.  In particular, for $I=(k,k,\dots,k)$
($r$ repetitions), we have
$$
A_I(p-1)\equiv S_I(p-1)\equiv 0\mod p
$$
for prime $p>rk+1$ (cf. Theorem 1.5 of \cite{Zh}).  There is
the following result relating sums associated to $I$ and $\bar I$
(cf. Lemma 3.2 of \cite{Zh}).
\begin{thm} For any composition $I$, $A_I(p-1)\equiv 
(-1)^{|I|}A_{\bar I}(p-1)\mod p$, and similarly
$S_I(p-1)\equiv(-1)^{|I|}S_{\bar I}(p-1)\mod p$.
\end{thm}
\begin{proof} Let $I=(i_1,\dots,i_k)$.  Working mod $p$, we have
\begin{multline*}
A_I(p-1)\equiv\sum_{p>a_1>\dots>a_k>0}\frac1{a_1^{i_1}\cdots a_k^{i_k}}
\equiv\sum_{p>a_1>\dots>a_k>0}\frac{(-1)^{i_1+\dots+i_k}}
{(p-a_1^{i_1})\cdots (p-a_k^{i_k})}\\
\equiv\sum_{0<b_1<\dots<b_k<p}\frac{(-1)^{i_1+\dots+i_k}}
{b_1^{i_1}\cdots b_k^{i_k}}=(-1)^{|I|}A_{\bar I}(p-1),
\end{multline*}
and similarly for $S_I$.
\end{proof}
An immediate consequence is that $S_I(p-1)\equiv A_I(p-1)\equiv 0\mod p$ if
$I=\bar I$ and $|I|$ is odd.  Another consequence is that
$S_{(i,j)}(p-1)\equiv A_{(i,j)}(p-1)\equiv 0\mod p$ when $p>i+j+1$
and $i+j$ is even.  This is because
$$
S_{(i,j)}(p-1)+S_{(j,i)}(p-1)\equiv 0\mod p
$$
for $p>i+j+1$ by Theorem 6.1, while $S_{(i,j)}(p-1)\equiv S_{(j,i)}(p-1)
\mod p$ when $i+j$ is even by Theorem 6.6.
\par
We have the following result relating $S_I$ and $S_{I^*}$.
\begin{thm} $S_I(p-1)\equiv -S_{I^*}(p-1)\mod p$ for all primes $p$.
\end{thm}
\begin{proof} Let $f$ be a sequence.  From the definition of $\del$
$$
\Si\del f(n)=\sum_{i=0}^n\binom{n+1}{i+1}(-1)^if(i),
$$
so taking $n=p-1$ gives
$$
\Si\del f(p-1)\equiv (-1)^{p-1}f(p-1)\equiv f(p-1)\mod p .
$$
Now take $f=S_I$ and apply Theorem 6.4.
\end{proof}
\par
This result has the following corollary for the $A_I$, which
may be compared with Theorem 4.4 of \cite{H1}.
(We use superscripts for repetition, so $(n,1^k)$ means
the composition of weight $n+k$ with $k$ repetitions of 1.)
\begin{thm} If $p$ is a prime with $p>\max\{k+1,n\}$, then
$$
A_{(n,1^k)}(p-1)\equiv A_{(k+1,1^{n-1})}(p-1)\mod p .
$$
\end{thm}
\begin{proof} First note that $(n,1^k)^*=(1^{n-1},k+1)$.  So, combining
Theorems 6.7 and 6.6,
\begin{equation}
S_{(n,1^k)}(p-1)\equiv -S_{(1^{n-1},k+1)}(p-1)
\equiv (-1)^{n+k+1}S_{(k+1,1^{n-1})}(p-1)\mod p .
\end{equation}
Now equate the right-hand sides of parts (1) and (2) of Theorem 6.2 and
then apply $\rho_{p-1}$ to get
$$
(-1)^{\ell(I)}S_{\bar I}(p-1)=\sum_{I_1\sqcup\dots\sqcup I_l=I}(-1)^l
A_{I_1}(p-1)\cdots A_{I_l}(p-1)
$$
for any composition $I$; if we set $I=(1^k,n)$, the hypothesis insures
that all the terms on the right-hand side are zero mod $p$ except the 
one with $l=1$, giving $(-1)^kS_{(n,1^k)}(p-1)\equiv A_{(1^k,n)}(p-1)\mod p$.  
Apply Theorem 6.6 to get $S_{(n,1^k)}(p-1)\equiv (-1)^n A_{(n,1^k)}(p-1)
\mod p$.  By the same argument, $S_{(k+1,1^{n-1})}(p-1)\equiv 
(-1)^{k+1}A_{(k+1,1^{n-1})}(p-1)\mod p$, and equation (9) gives the
conclusion.
\end{proof}
\par
We can state Theorem 6.7 in algebraic language as follows.
Define, for each prime $p$, a map $\chi_p:\H_{\mathbf Z}^1\to\Zp$ by
$\chi_p(w)=\rho_{p-1}(w)$.  (Here $\H_{\mathbf Z}^1$ is the
integral version of $\H^1$, i.e., the graded $\mathbf Z$-module in
$\Zxy$ generated by words ending in $y$.)  The commutative diagram
(8) gives the following algebraic version of Theorem 6.7, which can
be considered a mod $p$ counterpart of the duality theorem for MZVs.
\begin{thm} As elements of $\Zp$, $\chi_p(w)=\chi_p(\psi(w))$
for words $w$ of $\H^1$.
\end{thm}
For example, since $\psi(x^2y^3)=-x^2y^3-xy^4-yxy^3-y^5$, we have
$$
A_{(3,1,1)}(p-1)\equiv -A_{(3,1,1)}(p-1)-A_{(2,1,1,1)}(p-1)
-A_{(1,2,1,1)}(p-1)-A_{(1,1,1,1,1)}(p-1)\mod p .
$$
For $p>6$ this reduces to $2A_{(3,1,1)}(p-1)\equiv 
-A_{(2,1,1,1)}(p-1)-A_{(1,2,1,1)}(p-1)\mod p$.

\end{document}